\documentclass[11pt,oneside,a4paper]{article}
\usepackage{hyperref}
\usepackage{subcaption}
\usepackage{graphicx}
\usepackage{amsfonts}
\usepackage{graphicx}
\usepackage{epstopdf}
\usepackage{algorithmic}
\usepackage{mismath}
\usepackage{cancel}
\usepackage{amsmath}
\usepackage{amsthm}
\usepackage{float}
\theoremstyle{definition}
\newtheorem{Def}{Definition}

\newtheorem{remark}[Def]{Remark}
\newtheorem{theorem}[Def]{Theorem}

\newcommand{\U}{\mathbf{U}}
\newcommand{\Vee}{\mathbf{V}}

\newcommand{\demi}{\frac{1}{2}}

\newcommand{\yu}{\mathbf{u}}

\newcommand{\rhoti}{\tilde{\rho}}

\newcommand{\ut}{\tilde{u}}

\newcommand{\ENSK}{E_{\textnormal{NSK}}}
\newcommand{\Ld}{\Delta_h}
\newcommand{\Lp}{\Delta}

\DeclareMathOperator{\divn}{div}

\newcommand{\gradx}{\nabla}

\newcommand{\divx}{\divn}

\newcommand{\fx}{\nabla^+_x}
\newcommand{\fy}{\nabla^+_y}
\newcommand{\cx}{\nabla^c_x}
\newcommand{\cy}{\nabla^c_y}
\newcommand{\bx}{\nabla^-_x}
\newcommand{\by}{\nabla^-_y}

\newcommand{\tdt}{\frac{\di}{\di t}}

\ifpdf
\DeclareGraphicsExtensions{.eps,.pdf,.png,.jpg}
\else
\DeclareGraphicsExtensions{.eps}
\fi

\title{A Structure Preserving Finite Volume Scheme for the Navier-Stokes-Korteweg Equations\thanks{All autors were supported by the DFG within the priority research program SPP 2410, project 
		525866748. 
}}
\author{Jan Giesselmann\thanks{Department of Mathematics, Technische Universität Darmstadt, Germany ({giesselmann@mathematik.tu-darmstadt.de})}
	\and  Philipp \"Offner\thanks{Institute of Mathematics, Clausthal University of Technology, Germany, ({philipp.oeffner@tu-clausthal.de})}
	\and Robert Sauerborn\thanks{Institute of Mathematics, Clausthal University of Technology, Germany, ({robert.sauerborn@tu-clausthal.de})}
}

\begin{document}

\maketitle

\begin{abstract}
We present a semi-discrete finite volume scheme for the local Navier-Stokes-Korteweg and Euler-Korteweg systems. Our scheme is applicable for equidistant Cartesian meshes in one and two space dimensions. {In contrast to other works, which employ, for example, hyperbolic approximations of the equations or auxiliary-variable approaches leading to extended systems, our scheme operates directly on the original system.} We prove that it conserves mass and momentum and is energy stable. Numerical experiments complement our theoretical findings, showing that the scheme is convergent of order one if employed with explicit or implicit time discretisation.
\end{abstract}

\section{Introduction} 
\label{sec:introduction}

We study numerical schemes for the (local) isothermal Navier-Stokes-Korteweg (NSK) system
\begin{subequations}\label{eq_Korteweg}
	\begin{align}
		\partial_t \rho + \divx (\rho \yu) =&\, 0 \label{eq_Korteweg_rho}\\
		\partial_t (\rho\yu) + \divx\left[ \rho \yu \otimes \yu + p(\rho )\mathbf{I}\right] =&\, \mu \Delta \yu + \kappa \divx\left[  \left( \rho \Lp\rho + \demi|\gradx\rho|^2\right) \mathbf{I} - \gradx \rho \otimes \gradx \rho\right] \label{eq_Korteweg_m} \notag\\
		~
	\end{align}
\end{subequations}
in the conservative variables density, $\rho>0$, and momentum, $\rho \yu$, with a pressure function $\rho \mapsto p(\rho)$, and the simplified viscosity tensor $\nabla\yu$. We employ this simplification as the focus of this paper is on discretizing the Korteweg stress tensor. We consider a constant capillarity coefficient $\kappa(\rho) = \kappa>0$ and constant viscosity $\mu \geq 0$. The system \eqref{eq_Korteweg} is called the Euler-Korteweg (EK) system for $\mu=0$. We consider \eqref{eq_Korteweg_rho}, \eqref{eq_Korteweg_m} endowed with suitable initial data and periodic boundary conditions. 

This is a well known model for, e.g., compressible multi phase flows, see e.g. \cite{Rohde2010} and thin film shallow water flows, see e.g. \cite{noblevilaEKthinfilmflow}. It is also used in many applications for single phase flows and, arguably, captures more complex phenomena than the Navier-Stokes equations \cite{SlemrodReview,zbMATH00013316}. We focus on the case of a continuously differentiable and strictly monotone pressure law so that the first order part is strictly hyperbolic and any potential energy $P$ satisfying the differential equation $P'(\rho)\rho - P(\rho) = p(\rho)$ is convex. 

For this system several numerical schemes have already been developed. Let us start our discussion with those that provide energy consistency, i.e., the schemes dissipate energy and all dissipation comes from physical viscosity. The series of papers \cite{noblevilaEKthinfilmflow,zbMATH06579844,zbMATH07507231} derived mass and  momentum conserving and energy consistent schemes (on Cartesian meshes) which  are based on an augmented formulation that uses an evolution equation for $\sqrt{\kappa} \partial_x \rho / \sqrt{\rho}$. 
Recently, a mass and momentum conservative and energy dissipative scheme was introduced for a hyperbolic approximation to the NSK system \cite{zbMATH07599597, dhaoadi11040876}. Previously elliptic and parabolic relaxations of the third order terms had been considered \cite{keim2022relaxation,zbMATH07318386} which made the first order part of the problem hyperbolic, even for multiphase flows. Thus, structure preserving methods for hyperbolic problems could be applied. There are also a variety of schemes that are asymptotic preserving in the low Mach limit \cite{zbMATH06818400,MR3335225}.

In addition, there is a variety of schemes for the more demanding case of non-monotone pressure where the simultaneous conservation of mass and momentum as well as dissipation of energy has not been achieved yet. Let us refer to \cite{janenergy2014, zbMATH06004288, zbMATH08071323} for schemes that sacrifice conservation of momentum to ensure energy consistency even for non-monotonous pressure functions.

The main contribution of this paper is to present a numerical scheme that provides mass and momentum conservation and energy consistency without introducing an extended system. Our scheme is intended as a fallback scheme for structure preserving higher order schemes that rely on limiters \cite{kuzmin2024} or on the MOOD approach \cite{clain2011}.
Currently there are no convergence proofs for numerical schemes for the NSK system if existence of a smooth solution is not assumed a priori. We believe that our scheme is a good candidate for a first convergence proof in the spririt of dissipative measure-valued (DMV) solutions, see \cite{FeireislEtAl2021}. Such a proof would also establish existence of DMV solutions for the NSK and EK systems.

The remainder of this paper is organized as follows: In Section 2, we state the scheme in one space dimension and for Cartesian meshes in two space dimensions. In Section 3, we discuss mass and momentum conservation as well as energy consistency -- all in the case of periodic boundary conditions. The proofs for the two-dimensional case are postponed to the appendix. Finally, we provide numerical experiments in Section 4. These confirm the energy consistency and suggest that the scheme is first order convergent. We also study the behaviour of solutions in the $\kappa \rightarrow 0$ limit.

\section{Semidiscrete Methods}
\subsection{Semidiscretization of the 1D NSK system} 
\label{sec:1Dsemidis} 
In order to present our strategy in a technically simpler setting, we first consider the isothermal NSK equations in one space dimension given by
\begin{subequations}\label{eq_Korteweg1D}
	\begin{align}
		\partial_t \rho+ \partial_x  (\rho u) &=0  \\
		\partial_t (\rho u) + \partial_x\left[  \rho u^2 +p(\rho)  \right]  &=  \mu \partial_{xx} 
		u + \kappa \partial_x \left[ \rho \partial_{xx}\rho - \demi\left( \partial_x\rho\right) ^2\right].
	\end{align}
\end{subequations}
We consider the one-dimensional torus $\mathbf{T} = \R / \Z$ divided into $N$ equally sized cells as the computational domain and propose the following finite volume method with a Lax-Friedrichs numerical flux:
\begin{subequations}\label{semi1d}
	\begin{align}
		-\frac{\di}{\di t} \rho_{i}=&\,     \cx(\rho_i u_i) -\lambda h \Ld\rho_i  \label{semirho1}\\	 \label{semiu1}
		-\frac{\di}{\di t}\left(  \rho_i u_i \right) =& \,    \cx(\rho_i u_i^2) + \cx p(\rho_i) - \lambda h 	\Ld(\rho_i u_i) - \mu \Ld u_i \\
		&\, - \kappa\bx\left[ \frac{\rho_{i+1}\Ld\rho_{i} + \rho_{i}\Ld\rho_{i+1}}{2} - \demi (\fx\rho_{i})^2\right] \notag
	\end{align}
\end{subequations}
Here $\phi_i$ designates the finite volume approximation to the cell average of the variable $\phi$ in the $i$-th cell $\left[ \frac{i-1}{N},\frac{i}{N}\right]$. (Since no reconstruction step is performed, one might as well regard $\phi$ as a piecewise constant function.)\\
The operators are defined as follows:
\begin{align*}
	\cx\phi_{i} &\coloneqq \frac{\phi_{i+1} - \phi_{i-1}}{2h}; \quad \nabla^\pm_x\phi_{i} \coloneq \pm \frac{\phi_{i\pm 1} - \phi_{i}}{h};\\
	\Ld\phi_{i} &\coloneqq \frac{\phi_{i+1}-2\phi_{i}+\phi_{i-1}}{h^2};
\end{align*}
with the cell size $h \coloneqq N^{-1}$ and the periodic boundary conditions realized as $\phi_{N+1} \coloneqq \phi_1,\; \phi_{0} \coloneqq \phi_N,$ etc. The following identities hold for these operators:
\begin{subequations}\label{rules}
\begin{align}
	\cx\phi_{i} =& \demi \left( \fx \phi_{i} + \bx \phi_{i} \right),\\
	 \Ld\phi_{i} =& \fx \bx \phi_{i} = \bx \fx \phi_{i},\\
	 \sum_{i=1}^{N} \cx\phi_{i} =& \sum_{i=1}^{N} \nabla^\pm_x\phi_{i} = \sum_{i=1}^{N} \Ld\phi_{i} = 0 \txt{(conservation property),}\label{eq_con_prop}\\
	  \sum_{i=1}^{N} \phi_{i} \fx \psi_i =& -\sum_{i=1}^{N} \psi_{i} \bx \phi_i \txt{(partial summation).}\label{eq_part_sum}
\end{align}
\end{subequations}
 Initial conditions are realized by evaluation at cell centres: $\rho_i(0) = \rho_0(x_i), \, u_i(0) = u_0(x_i)$ with $x_i = \frac{2i-1}{2N}$. The numerical flux coefficient $\lambda$ depends on the maximum of propagation speeds across all cells; $\lambda = \demi \max_{i= 1, \ldots, N} \left\lbrace |u_i| + \sqrt{p'(\rho_i)}\right\rbrace$. 
We will abbreviate the right-hand sides of (\ref{semirho1}) and (\ref{semiu1}) as $F^\rho(\rho_{i}, \rho_{i} u_i)$ and $F^{\rho u}(\rho_{i}, \rho_{i}u_i)$ respectively where it is convenient.\\
We note the ``cross-average" between neighbouring cells used in the discretization of the third order capillarity term. This is to achieve the balance of energy that we calculate in Section \ref{sec:Properties1D}, but poses challenges in the analysis of convergence to dissipative measure valued solutions in the spirit of \cite{FeireislEtAl2021}. We will address these challenges and carry out the convergence analysis in a future paper.

\begin{remark}
	In Section \ref{sec:numtest}, we also perform numericals tests with a Rusanov flux, where the numerical diffusion is weighted with the local propagation speed instead of a global maximum. The diffusion term then reads $-\left(  \lambda_{i+1/2} \fx \phi_i - \lambda_{i-1/2} \bx \phi_{i}\right)  $, with $\lambda_{i+1/2} = \demi \max_{k = i, i+1} \left\lbrace |u_k| + \sqrt{p'(\rho_k)}\right\rbrace$, instead of $-\lambda h \Ld\phi_i$. 
\end{remark}

\subsection{Semidiscretization of the 2D NSK system}
For the semidiscretization of the NSK system \eqref{eq_Korteweg} in two space dimensions, with $\yu = \left[ u,v\right]^T$ we consider the two-dimensional torus $\mathbf{T}^2$ divided in square grid cells\footnote{This is regrettably necessary, as only square grids produce conveniently quadratic terms in calculating the energy balance, see Appendix \ref{sec:AppA}.} with width and height $h \coloneqq N^{-1}$ for some $N \in \N$ as our computational domain, and define the following operators for cell-averaged variables $\phi$ on it:

\begin{align*}
	\Ld(\phi_{i,j}) \coloneqq \Delta_{h,x}\phi_{i,j} + \Delta_{h,x}\phi_{i,j} &\coloneqq \frac{\phi_{i+1,j}-2\phi_{i,j}+\phi_{i-1,j}}{h^2}+\frac{\phi_{i,j+1}-2\phi_{i,j}+\phi_{i,j-1}}{h^2}\\
	\nabla^c\phi_{i,j} &\coloneqq		
	\begin{pmatrix}
		\cx\phi_{i,j}  \\
		\cy\phi_{i,j}
	\end{pmatrix} \coloneqq
	\begin{pmatrix}
		\frac{\phi_{i+1,j} - \phi_{i-1,j}}{2h}  \\
		\frac{\phi_{i,j+1} - \phi_{i,j-1}}{2h}
	\end{pmatrix}\\
	\nabla^\pm\phi_{i,j} &\coloneq \begin{pmatrix}
		\nabla^\pm_x\phi_{i,j}  \\
		\nabla^\pm_y\phi_{i,j}\end{pmatrix} \coloneqq \pm \begin{pmatrix}
		\frac{\phi_{i\pm 1,j} - \phi_{i,j}}{h}  \\
		\frac{\phi_{i,j \pm1} - \phi_{i,j}}{h}
	\end{pmatrix}
\end{align*}
Boundary conditions are realized as $\phi_{i+kN, j+lN} = \phi_{i,j}$ for $k,l \in \mathbf{N}$ and initial conditions are evaluated at cell centres, as in the one-dimensional method. Calculation rules analogous to those for the one-dimensional operators also apply.\\
The semidiscrete method for the two-dimensional NSK system is then given as
\begin{subequations}\label{semi2d}
	\begin{align}
		- \frac{\di}{\di t}  \rho_{i,j} =&\,  \cx(\rho u)_{i,j}+ \cy(\rho v)_{i,j}-\lambda h \Delta_{h}\rho_{i,j} \label{semi2drho}\\	
		-\frac{\di}{\di t}  (\rho u)_{i,j}   =&\, \cx(\rho u^2)_{i,j}+ \cy(\rho u v)_{i,j} + \cx p(\rho_{i,j})  - \mu \Ld u_{i,j}\notag \\
		&-\lambda h\Delta_{h}(\rho u)_{i,j}\notag \\ 
		&- \kappa \Bigg[ \bx \frac{\rho_{i,j}\Ld\rho_{i+1,j} + \rho_{i+1,j}\Ld\rho_{i,j}}{2}-\demi \bx \left( \fx\rho_{i,j}\right)^2\ \\
		&\quad + \demi \bx\left(  \by\rho_{i+1,j}\by \rho_{i,j}\right)   - \by \left( \cx\rho_{i,j}\, \fy \rho_{i,j}\right) \Bigg] \notag \\
		-\frac{\di}{\di t}  (\rho v)_{i,j}   =&\, \cy(\rho v^2)_{i,j}+ \cx(\rho u v)_{i,j} + \cy p(\rho_{i,j}) - \mu \Ld v_{i,j} \notag \\		
		&-\lambda h\Delta_{h}(\rho v)_{i,j}\notag \\ 
		&- \kappa \Bigg[ \by \frac{\rho_{i,j}\Ld\rho_{i,j+1} + \rho_{i,j+1}\Ld\rho_{i,j}}{2}-\demi \by \left( \fy\rho_{i,j}\right)^2 \label{semi2dv} \\
		&\quad + \demi \by\left(  \bx\rho_{i,j+1}\,\bx \rho_{i,j}\right)   - \bx \left( \cy\rho_{i,j}\, \fx \rho_{i,j}\right) \Bigg].  \notag
	\end{align}
\end{subequations}
In addition to cross-averaged third order terms, the two-dimensional method also employs a mixture of forward, backward and central differencing operators as well as index shifts in the discretization of capillarity terms, posing more challenges for convergence analysis. This is however the simplest configuration we considered for which we could prove energy dissipation.

\section{Properties of proposed schemes}
In this section we prove the structure preserving properties of the one-dimensional scheme. The corresponding properties of the two-dimensional scheme are stated and their proofs are provided in Appendix \ref{sec:AppA}.
\subsection{Properties of the 1D method}\label{sec:Properties1D}

\begin{theorem}\label{conservation1d}
	The method (\ref{semi1d}) conserves the total mass and momentum given by $\frac{1}{N} \sum_{i=1}^{N} \rho_{i}$ and $\frac{1}{N} \sum_{i=1}^{N} \rho_{i}u_i$ respectively.
\end{theorem}
\begin{proof}
	
	This is a direct consequence of the conservative formulation of the method. Since all terms in the right-hand sides of (\ref{semirho1}) and (\ref{semiu1}) are multiples of operators with the conservation property \eqref{eq_con_prop} applied to variables under periodic boundary conditions, we have 
	\begin{align*}
		\frac{\di}{\di t} \sum_{i=1}^{N} \rho_{i} = \frac{\di}{\di t} \sum_{i=1}^{N} \rho_i u_i  = 0.
	\end{align*}
\end{proof}
\begin{theorem}\label{1DDissipation}
	The method (\ref{semi1d}) with a numerical dissipation parameter  $\lambda = \demi \max_{i= 1, \ldots, N} \left\lbrace |u_i| + \sqrt{p'(\rho_i)}\right\rbrace$ dissipates the total energy of the discretized NSK system given by
	\begin{align}
		\ENSK[\rho,\rho u]= \frac{1}{N} \sum_{i=1}^{N} E_i \txt{with} E_i \coloneqq \demi \rho_i u_i^2 + P(\rho_i) + \frac{\kappa}{2} \left( \fx\rho_i\right) ^2.
	\end{align}
In particular, the inequality
\begin{align}
	\frac{\di}{\di t}\ENSK[\rho,\rho u] \leq - \mu\sum_{i=1}^{N}  \left(  \fx u_i\right)^2 - \kappa \sum_{i=1}^{N} \lambda h (\Ld\rho_i)^2 \leq 0
\end{align}holds.
\end{theorem}
\begin{proof}
	We multiply equations (\ref{semirho1}) and (\ref{semiu1}) by the discrete entropy variables with respect to $\ENSK$, i.e. $P'(\rho_i) - \demi u_i^2 - \kappa\Ld\rho_i$ and $u_i$ respectively. Summing both over all cells yields
	\begin{align}\label{1DEbalance}
		-N \frac{\di}{\di t} \ENSK(\rho,u) = \sum_{i=1}^{N} \Bigg[ &\left( \cx(\rho_i u_i) -\lambda h \Ld\rho_i\right) \left(P'(\rho_i) - \demi u_i^2\right)\\
	\notag	&+ u_i \left( \cx(\rho_i u_i^2) + \cx p(\rho_i) - \lambda h \Ld(\rho_i u_i) - \mu \Ld u_i\right) \\
	\notag	&+ \kappa \Big[\lambda h \left( \Ld\rho_i\right) ^2 -\Ld\rho_i\cx (\rho_i u_i) +\demi u_i \bx(\fx\rho_i)^2 \\
	\notag	&- u_i \bx \frac{\rho_{i+1}\Ld\rho_{i} + \rho_{i}\Ld\rho_{i+1}}{2} \Big] \Bigg],
	\end{align}
	where we use the identity \begin{align*}\frac{\di}{\di t} \sum_{i=1}^{N} \demi \left( \fx\rho_i\right) ^2 =&\, h^{-2} \sum_{i=1}^{N} \frac{\di \rho_{i+1}}{\di t} (\rho_{i+1}-\rho_{i}) - \frac{\di \rho_{i}}{\di t} (\rho_{i+1}-\rho_{i})\\
	=&\,h^{-2 }\sum_{i=1}^{N} \frac{\di \rho_{i}}{\di t}\left[ (\rho_i -\rho_{i-1})-(\rho_{i+1}-\rho_{i})\right] = \sum_{i=1}^{N} \frac{\di \rho_{i}}{\di t}\, \Ld \rho_i.\end{align*}
	We split the terms on the right-hand side of equation (\ref{1DEbalance}) into
	\begin{align*}
		A \coloneqq \sum_{i=1}^{N} \Bigg[ &\left( \cx(\rho_i u_i) -\lambda h \Ld\rho_i\right) \left(P'(\rho_i) - \demi u_i^2\right)\\
		\notag	&+ u_i \left( \cx(\rho_i u_i^2) + \cx p(\rho_i) - \lambda h \Ld(\rho_i u_i) \right)\Bigg],\\
		B \coloneqq  - \mu\sum_{i=1}^{N} & u_i  \Ld u_i\\
		C \coloneqq \sum_{i=1}^{N} \kappa \Bigg[&\lambda h \left( \Ld\rho_i\right) ^2 -\Ld\rho_i\cx (\rho_i u_i) +\demi u_i \bx(\fx\rho_i)^2 \\
		\notag	&- u_i \bx \frac{\rho_{i+1}\Ld\rho_{i} + \rho_{i}\Ld\rho_{i+1}}{2} \Bigg],
	\end{align*}
	separating the contributions of viscosity ($B$) and capillarity ($C$) to the energy balance. We claim that $A\geq 0$, $B\geq 0$, and $C\geq 0$ hold, and thus $\frac{\di}{\di t} \ENSK[\rho, \rho u] \leq 0$.\\
	
	\emph{Proof of }$A\geq 0:$\\
 In \cite{TadmorEntropyStability2003}, entropy stability is proven for semi-discrete central difference schemes with Lax-Friedrichs numerical dissipation of the form 
	\begin{align}
		-\frac{\di}{\di t}\U_i = \cx \mathbf{F}(\U_i) - \lambda h \Ld \U_i,
	\end{align}
	provided that for the diffusion coefficient the lower bound 
	$$\lambda \geq \demi \max_{\U_i, k} \left\lbrace \left| \lambda_k\left( \frac{\di \mathbf{F}}{\di \U}(\U_i)\right) \right|\right\rbrace  $$
	holds, where $\lambda_k$ are the eigenvalues of the flux's Jacobian. Then, for any suitable convex entropy function $E(\U)$ with associated entropy variables $\Vee_i = \frac{\di E}{\di \U}(\U_i)$, we have
	\begin{align*}
		\sum_{i=1}^N \frac{\di E}{\di t}(\U_i) = \sum_{i=1}^N \Vee_i \cdot \frac{\di}{\di t} \U_i \leq 0.
	\end{align*} 
	We now consider the isothermal Euler system in conservative variables, together with its canonical entropy, i.e.:
	\begin{align*}
		\U_i &= (\rho_{i}, \rho_{i} u_i)^T, \\
		\mathbf{F}(\U_i) &= (\rho_{i} u_i, \rho_{i} u_i^2 + p(\rho_{i}))^T, \\
		E(\U_i) &= \demi \rho_{i} u_i^2 + P(\rho_{i}).
	\end{align*}
	Straightforward calculations show that with these choices of $\U, \mathbf{F}$, and $E$ we have
	\begin{align*}
		\sum_{i=1}^N \Vee_i \cdot \frac{\di}{\di t} \U_i &= -A\\
		\demi \max_{k} \left\lbrace \left| \lambda_k\left( \frac{\di \mathbf{F}}{\di \U}(\U_i)\right) \right|\right\rbrace  &= \demi \left(  |u_i| + \sqrt{p'(\rho_i)}\right)
	\end{align*}
	and thus $A \geq 0$ under the assumed lower bound on $\lambda$.\\
	
	\emph{Proof of }$B\geq 0:$\\
	Summation by parts gives us 
	\begin{align*}
	B = - \mu\sum_{i=1}^{N}  u_i  \Ld u_i = 	- \mu\sum_{i=1}^{N}  u_i  \bx \fx u_i = \mu\sum_{i=1}^{N}  \left(  \fx u_i\right)^2 \geq 0.
	\end{align*}
	While the calculations for $A$ and $B$ are well known standard results, the novelty of our method is that capillarity is discretized in a way that ensures energy stability without using auxilliary variables.

	\emph{Proof of }$C\geq 0:$\\
	We have $C = \kappa \sum_{i=1}^{N} \lambda h (\Ld\rho_i)^2 \geq 0$, as the remaining terms comprising $C$ sum up to zero:
		\begin{align*}
		&\sum_{i=1}^{N} - \Ld\rho_i\cx(\rho_i u_i) +\demi u_i \bx(\fx\rho_i)^2 \\
		=& \sum_{i=1}^{N} \Ld\rho_i \frac{\rho_{i+1}u_{i+1}-\rho_{i-1}u_{i-1}}{2h} +  \frac{u_i}{2h}\left( \frac{(\rho_{i+1}-\rho_i)^2}{h^2}-\frac{(\rho_i-\rho_{i-1})^2}{h^2}\right)\notag\\
		=& \sum_{i=1}^{N} u_i\left[ \frac{-\rho_i\Ld\rho_{i-1} + \rho_i\Ld\rho_{i+1}}{2h} + \frac{\rho_{i+1} - \rho_{i-1}}{2h} \frac{\rho_{i+1}-2\rho_{i}+\rho_{i-1}}{h^2}\right] \notag\\
		=& \sum_{i=1}^{N} u_i \left[ \frac{-\rho_i\Ld\rho_{i-1} + \rho_i\Ld\rho_{i+1}}{2h} + \frac{\rho_{i+1}\Ld\rho_i - \rho_{i-1}\Ld\rho_i}{2h}\right] \notag\\
		=& \sum_{i=1}^{N} u_i \bx \frac{\rho_{i+1}\Ld\rho_{i} + \rho_{i}\Ld\rho_{i+1}}{2}. 
	\end{align*}
	Indeed the first line and the last line with reversed sign in the above calculation are exactly (up to the factor of $\kappa$) the non-quadratic terms in $C$.
\end{proof}

\subsection{Properties of the 2D method}
{For the sake of readability, we only state the main results on our two-dimensional scheme \eqref{semi2d} here. The corresponding proofs are provided in Appendix \ref{sec:AppA}.  }
\begin{theorem}\label{2Dconservation}
	The semidiscrete method (\ref{semi2d}) conserves total mass and momentum given by $\frac{1}{N^2} \sum_{i,j=1}^{N} \rho_{i,j}$ and $\frac{1}{N^2} \sum_{i,j=1}^{N} \rho_{i,j}\yu_{i,j}$ respectively.
\end{theorem}
	\begin{theorem}\label{2Ddissipation}
	The semidiscrete method (\ref{semi2d}) with the numerical dissipation parameter $\lambda = \demi \max_{i,j} \left\lbrace |\yu_{i,j}| + \sqrt{p'(\rho_{i,j})}\right\rbrace$ dissipates the total discrete energy given by 
	\begin{align}\label{energy2d}
		\ENSK[\rho,\rho \yu]= \frac{1}{N^2} \sum_{i,j=1}^{ N} E_{i,j}; \; E_{i,j} \coloneqq \demi \rho_i |\yu_{i,j}|^2 + P(\rho_{i,j}) + \frac{\kappa}{2} \left| \nabla^+\rho_{i,j}\right| ^2.
	\end{align}
	In particular, the inequality 
	\begin{align}
		\frac{\di}{\di t}\ENSK[\rho,\rho \yu] \leq  - \mu \sum_{i,j=1}^{N} \left| \nabla^+ u_{i,j}\right|^2 + \left| \nabla^+ v_{i,j}\right|^2 - \kappa \sum_{i,j=1}^{N} \lambda h (\Ld\rho_{i,j})^2\leq 0
	\end{align}holds.	 
\end{theorem}

\section{Numerical Tests} 
\label{sec:numtest} 
In this section, we present the results of numerical experiments. They confirm convergence of the expected order for the proposed schemes, energy dissipation, and the correct approximation of the Euler system behaviour in the limit of vanishing capillarity and viscosity. We also replicate a numerical experiment on thin film flow performed with a different energy-dissipative scheme for the Euler-Korteweg system in \cite{noblevilaEKthinfilmflow} and compare results. Throughout this chapter, the time step $\Delta t$ in simulations is chosen dependent on the grid size $h$ according to the CFL-like condition 
$$ 
\Delta t = \alpha \left( \frac{\lambda}{h} + \frac{\mu}{h^2} + \frac{\kappa}{h^3}\right)^{-1},$$
where $\lambda$ is defined as half of the maximum absolute value of the eigenvalues to the hyperbolic part of the flux's Jacobian across the grid (as above), $\alpha=0.7$ for explicit time-stepping, and $\alpha = 20$ for implicit time-stepping.
\subsection{Convergence Studies}
\label{sub:num_mms} 
\subsubsection{1D case} 
We perform a convergence study using the method of manufactured solutions.  \\
The method of manufactured solutions (cf. chapter 6.3 in \cite{Verification}) is employed in numerical testing and verification when exact solutions to the equation underlying the numerical method to be tested are unavailable or impractical. The discretized system is forced along a prescribed manufactured solution $\left[ \rhoti(x,t), \rhoti \tilde{u}(x,t)\right]$ by solving the system of equations

	\begin{subequations}
		\begin{align*}
				\rho^{(n+1)}_i =& \, \rho^{(n)}_i - \Delta t F^\rho(\rho^{(n+1)}_i, \rho_i u^{(n+1)}_i) \\
				&\, + \Delta t\left[  \partial_t \rhoti + \partial_x \left( \rhoti\tilde{u}\right) \right]\Big|_{x_i,t^{(n+1)}}  \\
				\rho_i u_i^{(n+1)} =& \, \rho_i u^{(n)}_i - \Delta t F^{\rho u}(\rho^{(n+1)}_i, \rho_i u^{(n+1)}_i)\\
				&\, + \Delta t \Bigg[ \partial_t\left(  \rhoti \ut\right)  + \partial_x\left[  \rhoti \ut^2 +p(\rhoti)  \right]  -  \mu \partial_{xx} \ut - \\ &\,\kappa \partial_x \left[ \rhoti \partial_{xx} \rhoti - \demi\left( \partial_x\rhoti\right) ^2\right]\Bigg]\Bigg|_{x_i,t^{(n+1)}}
			\end{align*}
	\end{subequations}
to evolve the system from $t^{(n)}$ to $t^{(n+1)} = t^{(n)}+\Delta t$ in implicit Euler time-stepping; with initial conditions given by $\rho^{(0)}_i = \rhoti(x_i, 0)$ and $\rho_i u_i^{(0)} = \rhoti \tilde{u}(x_i,0)$.
Table \ref{MMS1D} shows the errors of the numerical solution $[\rho, \rho u]$ relative to the manufactured solution $[\rhoti, \rhoti \tilde{u}]$ given by
\begin{align*}
	\rhoti(x,t) =&\, 1 + 0.5\cos(2\pi x+t)\\
	\rhoti\tilde{u} (x,t) =&\, 0.5\sin(2\pi x+t)\rhoti(x,t)
\end{align*}
at $T=0.2$, with $\kappa = 0.01, \, \mu = 0.01$. The result suggests convergence of the order $h+\Delta t$ of the proposed one-dimensional method (\ref{semi1d}). The same experiment was also carried out for $\kappa \in \left\lbrace 10^{-4}, 10^{-3}, 10^{-1}\right\rbrace $ with virtually identical results.

\begin{table}[h]
	\centering
	$\begin{array}{|c|c|c|c|c|c|c|}
		\hline 
		\text{\#(grid cells)} & \text{rel. } L^1\text{ error of $\rho$} & EOC &\text{ rel. }L^1\text{ error of $\rho u$} & EOC  \\
		\hline
		32  & 0.0648 &      & 0.1946 & \\ 
		64  & 0.0344 & 0.963 & 0.1034 & 0.974  \\
		128 & 0.0180 & 1.052 & 0.0543 & 1.058  \\
		256 & 0.0089 & 1.062 & 0.0268 & 1.067  \\
		512 & 0.0044 & 1.020 & 0.0133 & 1.020  \\
		1024 & 0.0022 & 1.004 & 0.0067 & 1.004 \\
		\hline
	\end{array}$
	\caption{Convergence study for the 1D method (\ref{semi1d}), implicit Euler time-stepping}\label{MMS1D}
\end{table}
The method of manufactured solutions was also employed with explicit Euler time-stepping, so that the updates read
 	\begin{subequations}
 	\begin{align*}
 		\rho^{(n+1)}_i =& \, \rho^{(n)}_i - \Delta t F^\rho(\rho^{(n)}_i, \rho_i u^{(n)}_i) \\
 		&\, + \Delta t\left[  \partial_t \rhoti + \partial_x\left(  \rhoti\tilde{u}\right) \right]\Big|_{x_i,t^{(n)}}  \\
 		\rho_i u_i^{(n+1)} =& \, \rho_i u^{(n)}_i - \Delta t F^u(\rho^{(n)}_i, \rho_i u^{(n)}_i)\\
 		&\, + \Delta t \Bigg[ \partial_t \left( \rhoti \ut\right)  + \partial_x\left[  \rhoti \ut^2 +p(\rhoti)  \right]  -  \mu \partial_{xx} \ut\\ &\, - \kappa \partial_x \left[ \rhoti \partial_{xx}\rhoti - \demi\left( \partial_x\rhoti\right) ^2\right]\Bigg]\Bigg|_{x_i,t^{(n)}}.
 	\end{align*}
 \end{subequations}
The results are shown in Table \ref{MMS1DExpl}.  {Here, we recognize again convergence of order $h+\Delta t$. In summary, both the explicit and implicit Euler time-discretizations together with the semidiscretization yield the expected order of convergence. }Both experiments were also carried out for $\kappa \in \left\lbrace 10^{-4}, 10^{-3}, 10^{-1}\right\rbrace $ and with Rusanov numerical dissipation, with virtually identical results.
\begin{table}[h]
	\centering
	$\begin{array}{|c|c|c|c|c|c|c|}
		\hline 
		\text{\#(grid cells)} & \text{rel. } L^1\text{ error of $\rho$} & EOC &\text{ rel. }L^1\text{ error of $\rho u$} & EOC  \\
		\hline
		32  & 0.0642 &      & 0.1921 & \\ 
		64  & 0.0332 & 0.950 & 0.0991 & 0.950  \\
		128 & 0.0168 & 1.003 & 0.0500 & 0.996  \\
		256 & 0.0086 & 0.946 & 0.0260 & 0.938  \\
		512 & 0.0044 & 0.977 & 0.0132 & 0.974  \\
		1024 & 0.0022 & 0.991 & 0.0066 & 0.990 \\
		\hline
	\end{array}$
	\caption{Convergence study for the 1D method (\ref{semi1d}), explicit Euler time-stepping}\label{MMS1DExpl}
\end{table}
\subsubsection{2D case}
The proposed two-dimensional method (\ref{semi2d}) was also tested with the method of manufactured solutions and explicit Euler time-stepping. Table \ref{MMS2D} shows the resulting errors relative to the  manufactured solution $(\rhoti, \rhoti \tilde{u}, \rhoti \tilde{v})$ given by
\begin{align*}
	\rhoti(x,t) =&\, 0.5 + \sin^2(x+t) + \cos^2(y+t)\\
	\tilde{u} (x,t) =&\, \sin(x+t) \cos(y+t) \\
	\tilde{v} (x,t) =&\, \cos(x+t) \sin(y+t)
\end{align*}
at $T=0.2$, with $\kappa = 0.01, \, \mu = 0.01$, which again suggest convergence of order $h+\Delta t$.
\begin{table}[h]
	\centering
	$\begin{array}{|c|c|c|c|c|c|c|c|c|}
		\hline 
		\text{\#(grid cells)} & \text{rel. } L^1\text{ error of $\rho$} & EOC &\text{$\hdots$ of $\rho u$} & EOC &\text{$\hdots$ of $\rho v$} & EOC \\
		\hline
		32^2  & 0.03841 &      & 0.08052 & &0.10787&\\ 
		64^2  & 0.01986 & 0.952 & 0.04151 & 0.956  &0.05765&0.904\\
		128^2 & 0.00994 & 0.999 & 0.02086 & 0.993  &0.02933&0.975\\
		256^2 & 0.00514 & 0.951 & 0.01077 & 0.954  &0.01517&0.951\\
		512^2 & 0.00262 & 0.970 & 0.00555 & 0.956  &0.00771&0.976\\
		
		\hline
	\end{array}$
	\caption{Convergence study for the 2D method (\ref{semi2d}), explicit Euler time-stepping}\label{MMS2D}
\end{table}
{Virtually identical results have been obtained using explicit Euler time-stepping with Rusanov numerical dissipation.}

\subsection{Riemann Problem} 
\label{sub:num_riemann} 
\begin{figure}[H]
	\centering
%
%
		\includegraphics[width = 0.7\textwidth]{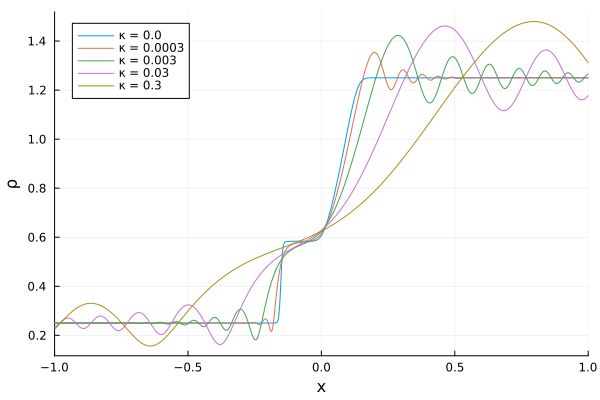}\caption{Riemann problem with contact discontinuity, densities at $t = 0.1$}\label{fig:Riemann1D020}

	\end{figure}

\begin{figure}[H]
		\centering
		\begin{subfigure}{0.495\textwidth}
		\includegraphics[width=\textwidth]{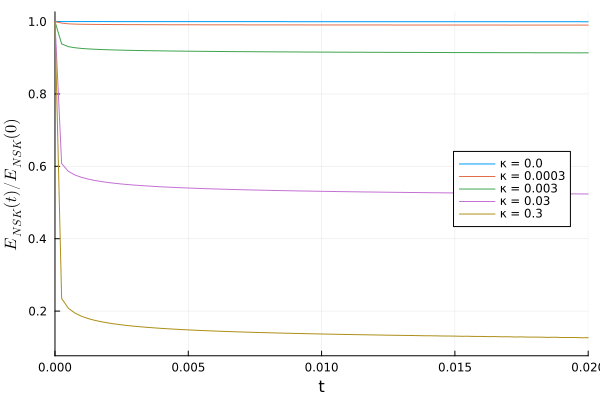}\caption{Evolution of energies}\label{fig:RiemannEnergy}

	\end{subfigure}\begin{subfigure}{0.495\textwidth}
	\includegraphics[width=\textwidth]{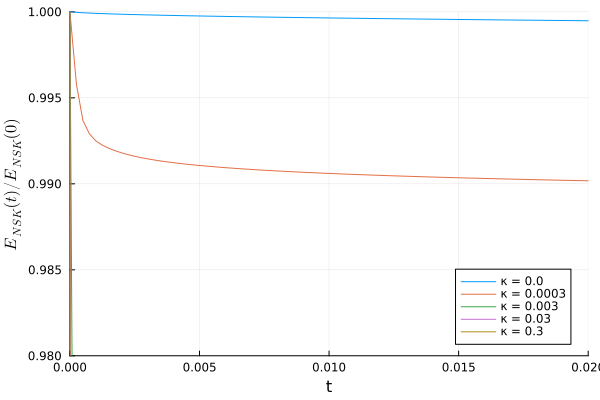}\caption{Evolution of energies (magnified)}\label{fig:RiemannEnergyZoom}

	\end{subfigure}

	\caption{Riemann problem with contact discontinuity}
\end{figure}

We tested the one-dimensional method (\ref{semi1d}) on the contact discontinuity Riemann problem $\rho_l = 0.25, \rho_r = 1.25, u_l = u_r = 0$ for the Euler-Korteweg system with explicit Euler time-stepping and Lax-Friedrichs numerical dissipation on a grid of size $h=2^{-10}$. 
The numerical density at $t = 0.1$ is shown in Figure \ref{fig:Riemann1D020}; with capillarity coefficients $\kappa$ of varying magnitude. Numerical solutions for low $\kappa$ closely resemble the barotropic Euler shock profile, however the shock is smeared out at higher capillarity coefficients. This is due to higher velocities being produced from initial capillary forces. The results seem to qualitatively suggest convergence of the approximations in the $L^1$-norm to a barotropic Euler shock profile for $\kappa, h \to 0$.\\
In Figures \ref{fig:RiemannEnergy} and \ref{fig:RiemannEnergyZoom}, the time evolution of discrete energies (denoted by $\ENSK$, as defined above) of numerical solutions relative to the initial energy are shown, with the upper left corner magnified to allow discerning the graphs for $\kappa = 0$ and $\kappa = 0.0003$. The increasing magnitude of surface tension's contribution to the initial discrete energy are evident, and while we have not proven dissipation of energy for the fully discretized method, the results indicate it is indeed dissipative. 

%
%
%
%

\subsection{Thin Film Experiment} 
\label{sub:num_sw} 
\begin{figure}[H]
	\centering
	\includegraphics[width=0.7\textwidth]{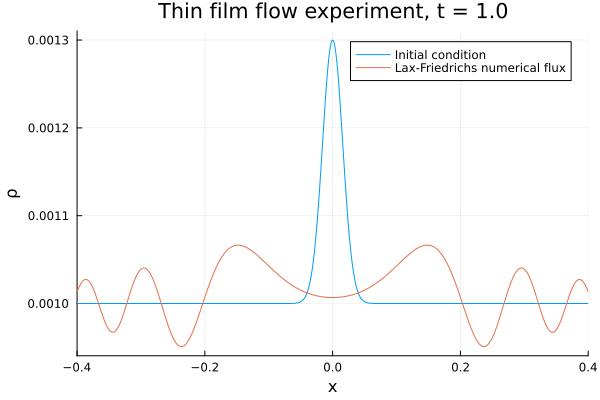}\caption{Profile of fluid surface at t = 1}\label{ThinFilm}
\end{figure}

We replicated the numerical experiment on thin film flow performed in \cite{noblevilaEKthinfilmflow} on a grid of size $h = 10^{-3}$, simulating a fluid on a horizontal plane without friction at the bottom. We set $\mu = 0, \, \kappa = 0.0059, \, p(\rho) = \frac{9.81}{2} \rho^2$, corresponding to a shallow water system with the capillarity coefficient of an aqueous solution of glycerin. Initial conditions are given by $\rho_0(x) = 10^{-3}(1+\exp(-2000x^2)), \, u_0(x) = 0$. The result using a Lax-Friedrichs numerical flux with explicit Euler time-stepping is shown in Figure \ref{ThinFilm}. The experiment was also performed using a Rusanov numerical flux with virtually indistinguishable results. \\ 
We note that while the positions of the central fluid peaks coincide with those in the simulation by Noble and Vila, our simulation results show stronger and faster propagating oscillations beyond the central fluid peaks. We speculate that this is because their auxiliary variable approach with three transport equations, each with their own numerical viscosity, might have a stronger dampening effect on oscillations. Increasing numerical dissipation in our method indeed generates results that resemble those in \cite{noblevilaEKthinfilmflow} more closely.
\newpage

\appendix 
\section{Proofs of structure preserving properties for the 2D method} \label{sec:AppA}
	\begin{proof}[Proof of Theorem \ref{2Dconservation}]
		Analogous to the proof of \ref{conservation1d}, the claim 
		\begin{align*}
			\frac{\di}{\di t} \sum_{i,j=1}^{N} \rho_{i,j} = \frac{\di}{\di t} \sum_{i,j=1}^{N} \left( \rho u\right)_{i,j} = \frac{\di}{\di t} \sum_{i,j=1}^{N} \left( \rho v\right)_{i,j} = 0
		\end{align*}
		
		follows from the conservative formulation of the method \eqref{semi2d}.
	\end{proof}

\begin{proof}[Proof of Theorem \ref{2Ddissipation}]
	Analogous to the proof of \ref{1DDissipation}, we multiply (\ref{semi2drho})-(\ref{semi2dv}) by $P'(\rho_{i,j}) -\demi |\yu_{i,j}|^2 - \Ld\rho_{i,j}, \; v_{i,j}$ and $u_{i,j}$ respectively and sum over all cells yielding
	\begin{align*}
	&\quad -N^2 \tdt \ENSK[\rho,\rho \yu] \notag \\
	&=  \sum_{i,j=1}^{N} \Bigg[\left( P'(\rho_{i,j})-\demi |\yu_{i,j}|^2\right) \left( \cx(\rho u)_{i,j} + \cy(\rho v)_{i,j} -\lambda h\Delta_{h}\rho_{i,j}\right) \notag \\ 
	&+u_{i,j} \big[ \cx\left[ (\rho u^2)_{i,j}+ p(\rho_{i,j})\right] +\cy(\rho u v)_{i,j}- \lambda h\Delta_{h}(\rho u)_{i,j} - \mu \Ld u_{i,j} \big] \\
		&+v_{i,j} \big[ \cy\left[ (\rho v^2)_{i,j}+ p(\rho_{i,j})\right]+ \cx(\rho u v)_{i,j}  -\lambda h\Delta_{h}(\rho v)_{i,j} - \mu \Ld v_{i,j} \big]  \notag \\
	& + \lambda \kappa h(\Delta_{h}\rho_{i,j})^2 - \kappa \Ld(\rho_{i,j}) \left( \cx(\rho u)_{i,j} + \cy(\rho v)_{i,j}\right)   \notag \\
	&- \kappa u_{i,j} \Bigg[ \bx \frac{\rho_{i,j}\Ld\rho_{i+1,j} + \rho_{i+1,j}\Ld\rho_{i,j}}{2}-\demi \bx \left( \fx\rho_{i,j}\right)^2\ \notag \\
	&\quad + \demi \bx\left(  \by\rho_{i+1,j}\by \rho_{i,j}\right)   - \by \left( \cx\rho_{i,j}\, \fy \rho_{i,j}\right)\Bigg] \notag \\
	&- \kappa v_{i,j}  \Bigg[ \by \frac{\rho_{i,j}\Ld\rho_{i,j+1} + \rho_{i,j+1}\Ld\rho_{i,j}}{2}-\demi \by \left( \fy\rho_{i,j}\right)^2\ \notag \\
	&\quad + \demi \by\left(  \bx\rho_{i,j+1}\,\bx \rho_{i,j}\right)   - \bx \left( \cy\rho_{i,j}\, \fx \rho_{i,j}\right)\Bigg]\Bigg]. \notag
	\end{align*}
	We split the right-hand side of that energy balance as follows:
	\begin{align*}
		A_1 \coloneqq& \sum_{i,j=1}^{N} \Bigg[\left( P'(\rho_{i,j})-\demi u_{i,j}^2\right)\left( \cx(\rho u)_{i,j} - \lambda h\Delta_{h,x}\rho_{i,j}\right) \notag \\ 
		& \qquad +u_{i,j} \left( \cx(\rho u^2)_{i,j} + \cx p(\rho_{i,j})-\lambda h \Delta_{h,x}(\rho u)_{i,j} \right)  \notag \Bigg], \\
		A_2 \coloneqq& \sum_{i,j=1}^{N} \Bigg[\left( P'(\rho_{i,j})-\demi v_{i,j}^2\right)\left( \cy(\rho v)_{i,j} -\lambda h\Delta_{h,y}\rho_{i,j}\right) \notag \\ 
		& \qquad +v_{i,j} \left( \cy(\rho v^2)_{i,j} + \cy p(\rho_{i,j})-\lambda h \Delta_{h,y}(\rho v)_{i,j} \right)  \notag \Bigg],
		\end{align*}		
		$$\begin{aligned}
		A_3 \coloneqq& \sum_{i,j=1}^{N} u_{i,j} \left( \cy(\rho u v)_{i,j} -\lambda h\Delta_{h,y}(\rho u)_{i,j} \right) -\demi u_{i,j}^2 \left[ \cy(\rho v)_{i,j} -\lambda h \Delta_{h,y}\rho _{i,j}\right],    \notag \\ 
		A_4 \coloneqq & \sum_{i,j=1}^{N} v_{i,j} \left( \cx(\rho u v)_{i,j} -\lambda h\Delta_{h,x}(\rho v)_{i,j} \right)  - \demi v_{i,j}^2\left[  \cx(\rho u)_{i,j} -\lambda h\Delta_{h,x}\rho_{i,j}\right] ,  \notag \\
		B\coloneqq&  -  \mu\sum_{i,j=1}^{N} u_{i,j} \Ld u_{i,j} +  v_{i,j} \Ld v_{i,j}, \\
		C \coloneqq& \sum_{i,j=1}^{N} \lambda \kappa h (\Ld\rho_{i,j})^2 ,\\
		D_1 \coloneqq& \, \kappa \sum_{i,j=1}^{N}  \Bigg[- \Ld(\rho_{i,j}) \cx(\rho u)_{i,j}  - u_{i,j} \Bigg[ \bx \frac{\rho_{i,j}\Ld\rho_{i+1,j} + \rho_{i+1,j}\Ld\rho_{i,j}}{2} \notag \\
		&\qquad \quad -\demi \bx \left( \fx\rho_{i,j}\right)^2 + \demi \bx\left(  \by\rho_{i+1,j}\by \rho_{i,j}\right)   - \by \left( \cx\rho_{i,j}\, \fy \rho_{i,j}\right)\Bigg]\Bigg],\\
		D_2 \coloneqq& \, \kappa \sum_{i,j=1}^{N}  \Bigg[- \Ld(\rho_{i,j}) \cy(\rho v)_{i,j} -v_{i,j}  \Bigg[ \by \frac{\rho_{i,j}\Ld\rho_{i,j+1} + \rho_{i,j+1}\Ld\rho_{i,j}}{2}\ \notag \\
		&\qquad \quad-\demi \by \left( \fy\rho_{i,j}\right)^2 + \demi \by\left(  \bx\rho_{i,j+1}\,\bx \rho_{i,j}\right)   - \bx \left( \cy\rho_{i,j}\, \fx \rho_{i,j}\right)\Bigg]\Bigg].
	\end{aligned}$$
	We obviously have $C \geq 0$, and we claim that $A_1, A_2, A_3, A_4, B, \geq 0$ and $D_1 = D_2 = 0$ also hold, thus $\frac{\di}{\di t}\ENSK[\rho,\rho \yu] \leq 0$.\\
	
	\emph{Proof of }$A_1, A_2 \geq 0:$\\
As in the proof of Theorem \ref{conservation1d}, these terms coincide with entropy balances of 1D barotropic Euler systems, so that each sum over $i$ or $j$ respectively will be positive given the assumption on $\lambda$.\\
	
	\emph{Proof of }$A_3, A_4 \geq 0:$\\	
We first consider $A_3$; $A_4$ can be treated analogously. For the terms approximating the physical flux we have
 \begin{align}
 	S_{i,j} \coloneqq &\, u_{i,j} \cy(\rho u v)_{i,j} - \demi u_{i,j}^2 \cy(\rho v)_{i,j} \notag\\
 	=&\, \frac{(\rho v )_{i,j+1}}{2h} \left(u_{i,j} u_{i,j+1} - \demi u_{i,j}^2\right) - \frac{(\rho v )_{i,j-1}}{2h} \left(u_{i,j} u_{i,j-1} - \demi u_{i,j}^2\right)\notag\\
 	=&\, \frac{(\rho v )_{i,j+1}}{4h} \left(  u_{i,j+1}^2 -  \left( u_{i,j+1} -u_{i,j}\right)^2 \right) - \frac{(\rho v )_{i,j-1}}{4h} \left(   u_{i,j-1}^2 -  \left( u_{i,j-1} -u_{i,j}\right)^2 \right) \notag\\
 	=&\, \demi \cy(\rho u^2 v)_{i,j} - \frac{1}{4h} \left(  (\rho v )_{i,j+1} (u_{i,j+1} - u_{i,j} )^2 - (\rho v )_{i,j-1}(u_{i,j} - u_{i,j-1})^2 \right).\label{mixedterms1}
 \end{align}
	For the numerical dissipation terms we have 
	\begin{align}
		T_{i,j} \coloneqq & \, u_{i,j}\Delta_{h,y}(\rho u)_{i,j} - \demi u_{i,j}^2 \Delta_{h,y}\rho_{i,j} \notag\\
		=& \, \frac{u_{i,j}}{h^2} ((\rho u )_{i,j+1} - 2(\rho u )_{i,j} + (\rho u )_{i,j-1}) - \frac{u_{i,j}^2}{2h^2} (\rho_{i,j+1} - 2\rho_{i,j}+ \rho_{i,j-1})\notag\\
		=& \, \frac{1}{2h^2} \big( (\rho u^2 )_{i,j+1} - 2(\rho u^2 )_{i,j} + (\rho u^2 )_{i,j-1}  \notag\\
	&	\quad\quad - \rho_{i,j+1}(u_{i,j+1} - u_{i,j})^2 - \rho_{i,j-1}(u_{i,j} - u_{i,j-1})^2\big)\notag\\
	=& \, \demi \Delta_{h,y} (\rho u^2 )_{i,j} - \frac{1}{2h^2}\left( \rho_{i,j+1}(u_{i,j+1} - u_{i,j})^2 + \rho_{i,j-1}(u_{i,j} - u_{i,j-1})^2\right).\label{mixedterms2}
	\end{align}

	Equations \eqref{mixedterms1} and \eqref{mixedterms2} together give us, after cancellation of the pure difference terms:
	\begin{align*}
	A_3  =& \sum_{i,j=1}^N S_{i,j} - \lambda h T_{i,j}\\
	=& \sum_{i,j=1}^N \frac{\lambda}{2 h}\left( \rho_{i,j+1}(u_{i,j+1} - u_{i,j})^2 + \rho_{i,j-1}(u_{i,j} - u_{i,j-1})^2\right)\\
	& \qquad- \frac{1}{4h} \left(  (\rho v )_{i,j+1} (u_{i,j+1} - u_{i,j} )^2 - (\rho v )_{i,j-1}(u_{i,j} - u_{i,j-1})^2 \right)\\
	\geq & \, \frac{1}{2 h} \sum_{i,j=1}^N \rho_{i,j+1}\left( \lambda  - \demi|v_{i,j+1}|\right) (u_{i,j+1} - u_{i,j} )^2\\
	& \quad \qquad +\rho_{i,j-1} \left( \lambda  - \demi|v_{i,j-1}|\right) (u_{i,j} - u_{i,j-1})^2\\
	\geq& \, 0, \txt{under the assumption $\lambda = \demi \max_{i,j} \left\lbrace |\yu_{i,j}| + \sqrt{p'(\rho_{i,j})}\right\rbrace$.}
	\end{align*}

	\emph{Proof of }$B\geq 0:$\\
	As in the one-dimensional case, summation by parts leaves us with sums over quadratic terms:
	$$ \begin{aligned}
		-  \mu\sum_{i,j=1}^{N} u_{i,j} \Ld u_{i,j} =& -  \mu\sum_{i,j=1}^{N} u_{i,j} \left( \bx\fx u_{i,j} + \by\fy u_{i,j}\right) \\ =&\, \mu \sum_{i,j=1}^{N} \left( \fx u_{i,j}\right) ^2 + \left( \fy u_{i,j}\right) ^2 \geq 0.
	\end{aligned}$$
	The calculation for the terms containing $v$ is analogous.\\

	
	\emph{Proof of }$D_1 = D_2 = 0:$\\
	We first consider $D_1$, $D_2$ can be treated analogously. We calculate a cancellation of terms containing exclusively differences in the $x$-direction as we did in the one-dimensional case: 	
	
	\begin{align*}
		&\sum_{i,j=1}^{N}  - \Ld\rho_{i,j} \cx(\rho u)_{i,j}  + \demi u_{i,j} \bx \left( \fx\rho_{i,j}\right)^2  \notag \\
		=& \sum_{i,j=1}^{N} u_{i,j} \left( \frac{-\rho_{i,j}\Ld\rho_{i-1,j} + \rho_{i,j}\Ld\rho_{i+1,j}}{2h_x} + \frac{\rho_{i+1,j} - \rho_{i-1,j}}{2h_x} \frac{\rho_{i+1,j}-2\rho_{i,j}+\rho_{i-1,j}}{h_x^2}\right) \notag\\
		=&\sum_{i,j=1}^{N} u_{i,j} \left( \bx \frac{\rho_{i,j}\Ld\rho_{i+1,j} + \rho_{i+1,j}\Ld\rho_{i,j}}{2} - \cx\rho_{i,j}\Delta_{h,y}\rho_{i,j}\right).
	\end{align*}
	The remaining terms in $D_1$ with differences in $x$- and $y$-direction finally cancel with $-\cx\rho_{i,j}\Delta_{h,y}\rho_{i,j}$:
	\begin{align}
		&\by \left(\cx\rho_{i,j}\, \fy \rho_{i,j}\right) - \demi \bx\left(  \by\rho_{i+1,j}\, \by \rho_{i,j}\right) \notag \\
		=& \,\frac{1}{2h_x h_y^2}\Big[\rho_{i+1,j}\left( \rho_{i,j+1} - \rho_{i,j}\right) - \rho_{i-1,j}\left( \rho_{i,j+1} - \rho_{i,j}\right) - \cancel{\rho_{i+1,j-1} \left( \rho_{i,j} - \rho_{i,j-1}\right)}\notag \\& \qquad+ \cancel{\rho_{i-1,j-1} \left( \rho_{i,j} - \rho_{i,j-1}\right)} \notag 
		-\rho_{i+1,j}\left(\rho_{i,j} - \rho_{i,j-1}\right) + \cancel{\rho_{i+1,j-1} \left( \rho_{i,j} - \rho_{i,j-1}\right)} \\
		&\qquad+ \rho_{i,j} \left( \rho_{i-1,j} - \cancel{\rho_{i-1,j-1}}\right) - \rho_{i,j-1}\left( \rho_{i-1,j} - \cancel{\rho_{i-1,j-1}}\right)\Big]\notag \\
		=&\, \frac{1}{2h_x h_y^2}\left[ \left( \rho_{i+1,j} - \rho_{i-1,j}\right) \left( \rho_{i,j+1} - 2\rho_{i,j} + \rho_{i,j-1}\right) \right] \notag \\
		=&\, \cx\rho_{i,j} \Delta_{h,y}\rho_{i,j}\notag.
	\end{align}
\end{proof}
 
\bibliographystyle{abbrv}
\bibliography{literature}
\end{document}